\documentclass[11pt]{article}
\usepackage{color,latexsym,amsmath,amssymb,path,url}

\newcommand{\ignore}[1]{}

\newtheorem{theorem}{Theorem}[section]
\newtheorem{lemma}[theorem]{Lemma}
\newtheorem{corollary}[theorem]{Corollary}

\newtheorem{conjecture}[theorem]{Conjecture}

\newcommand{\Proof}[1]
        {
        \noindent
        \emph{Proof #1.}~
        }
\newsavebox{\smallProofsym}                     
\savebox{\smallProofsym}                        %
        {
        \begin{picture}(7,7)                    %
        \put(0,0){\framebox(6,6){}}             %
        \put(0,2){\framebox(4,4){}}             %
        \end{picture}                           %
        }                                       %
\newcommand{\smalleop}[1]
        {
        \mbox{} \hfill #1~~\usebox{\smallProofsym}\!\!\!\!\!\!\
        }
\newenvironment{theProof}[1]
        {
        \Proof{#1}}{\smalleop{}
        \medskip
        }
\usepackage{graphicx,psfrag}

\newcommand{\RP}{\mathbb{R}\mathbf{P}}

\newcommand{\RR}{\ensuremath{\mathbb R}}
\newcommand{\CC}{\ensuremath{\mathbb C}}

\newcommand{\pts}{\mathcal P}
\newcommand{\ptss}{\mathcal V}
\newcommand{\curves}{\Gamma}

\begin{document}
\pagenumbering{arabic}

\title{Few distinct distances implies no heavy lines or circles
\thanks{
The first author was partially
supported by Grant 338/09 from the Israel Science Fund
and by the Israeli Centers of Research Excellence (I-CORE)
program (Center No.~4/11). Work by the second author was partially supported by the Department of Defense through the National Defense Science \& Engineering Graduate Fellowship (NDSEG) Program.
Work by the third author was partially supported by Swiss National Science Foundation Grant no. 200021-513047.
}
}


\author{
Adam Sheffer\thanks{%
School of Computer Science, Tel Aviv University,
Tel Aviv 69978, Israel.
{\sl adamsh@gmail.com}
}
\and
Joshua Zahl\thanks{%
Department of Mathematics, MIT, Cambridge, MA,  02139.
{\sl jzahl@mit.edu}
}
\and
Frank de Zeeuw\thanks{
EPFL, Lausanne, Switzerland.
{\sl fdezeeuw@gmail.com}
}
}

\maketitle

\begin{abstract}
We study the structure of planar point sets that determine a small number of distinct distances.
Specifically, we show that if a set $\pts$ of $n$ points determines $o(n)$ distinct distances, then no line contains $\Omega(n^{7/8})$ points of $\pts$ and no circle contains $\Omega(n^{5/6})$ points of $\pts$.

We rely on the bipartite and partial variant of the Elekes-Sharir framework that was presented by Sharir, Sheffer, and Solymosi in \cite{SSS13}. For the case of lines we combine this framework with a theorem from additive combinatorics, and for the case of circles we combine it with some basic algebraic geometry and a recent incidence bound for plane algebraic curves by Wang, Yang, and Zhang \cite{WYZ13}.
A significant difference between our approach and that of \cite{SSS13} (and other recent extensions) is that, instead of dealing with distances between two point sets that are restricted to one-dimensional curves, we consider distances between one set that is restricted to a curve and one set with no restrictions on it.
\end{abstract}

\section{Introduction}

Given a set $\pts$ of points in $\RR^2$, let $D(\pts)$ denote the number of distinct distances that are determined by pairs of points from $\pts$, and let $D(n)$ be the minimum of $D(\pts)$ over all $\pts\subset \RR^2$ with $|\pts|=n$.
In his celebrated 1946 paper \cite{erd46}, Erd\H os derived the upper bound $D(n) = O(n/\sqrt{\log n})$.
Recently, after 65 years and a series of increasingly larger lower bounds, Guth and Katz \cite{GK11} provided an almost matching lower bound $D(n) = \Omega(n/\log n)$.
In the process Guth and Katz developed several novel techniques, relying on tools from algebraic geometry.

The common belief is that $D(n) = \Theta(n/\sqrt{\log n})$.  While this leaves only a multiplicative gap of $\Theta(\sqrt{\log n})$, hardly anything is known about those sets that determine such a small number of distinct distances. To quote Gy\"orgy Elekes \cite{Elek02}:
\begin{quotation}\em ``Everywhere in mathematics, whenever we determine the maximum or minimum of a quantity, it is interesting to describe those configurations for which this extremum is attained. Sometimes even the stability of the extremal structures is interesting i.e., if we are not far from the best possible value, can the structure change very much or must it remain close to the optimal one? Such questions, of course, tend to be harder than those of the first type."
\end{quotation}
Erd\H os \cite{erd46} proved his upper bound by considering a $\sqrt{n}\times\sqrt{n}$ integer lattice, and the same bound is also attained by other types of $\sqrt{n}\times\sqrt{n}$ lattices.
It is shown in \cite{CSS13} that even a slightly uneven $n^{\alpha}\times n^{1-\alpha}$ integer lattice  (where $0 < \alpha < 1/2$) determines $\Theta(n)$ distinct distances, so such lattice structures appear to be somewhat unstable in the sense of Elekes's quote.
It remains unknown whether every set with $O(n/\sqrt{\log n})$ distinct distances ``has lattice structure", as Erd\H os conjectured \cite{erd86}.

As a first step in this direction, Erd\H os asked whether every such set can be covered by a small number of lines.
Since this also turned out to be difficult, Erd\H os suggested to prove that an optimal point set must have $\Omega(n^{1/2})$ points on a line, or even only $\Omega(n^\varepsilon)$ points.
Embarrassingly, after over 25 years, this problem is still wide open.
The only ``crack in the wall" was made by Szemer\'edi, who proved that for every set $\pts$ of $n$ points that determines $o(n)$ distinct distances there exists a line that contains $\Omega(\sqrt{\log n})$ points of $\pts$ (Szemer\'edi's proof was communicated by Erd\H os in \cite{erd75}; a variant of this proof can be found in \cite{PA95}, Theorem 13.7).

In this paper we make a crack on the other side of the wall, by proving that any set that comes close to being optimal cannot have too many points on a line.

\begin{theorem} \label{th:lines}
Let $\pts$ be a set of $n$ points in $\RR^2$ such that $D(\pts) = o(n)$. 
Then there does not exist a line containing $\Omega(n^{7/8})$ points of $\pts$.
\end{theorem}

The theorem cannot be extended to say anything about point sets $\pts$ for which $D(\pts) = o(n^\beta)$ for some $\beta>1$, since a line that contains $n^{\alpha}$ evenly spaced points (for any $7/8 \le \alpha < 1$) can be completed to an $n^{\alpha}\times n^{1-\alpha}$ integer lattice. As shown in \cite{CSS13}, such a lattice determines $\Theta(n)$ distinct distances.
We do believe that the constant $7/8$ in Theorem \ref{th:lines} can be further improved, and conjecture the following.

%
\begin{conjecture} \label{con:lines}
Let $\pts$ be a set of $n$ points in $\RR^2$ such that $D(\pts)=O(n/\sqrt{\log n})$.
Then there does not exist a line containing $\Omega(n^{1/2+\varepsilon})$ points of $\pts$, for any $\varepsilon>0$.
\end{conjecture}

A related ``bipartite'' problem is Purdy's conjecture, which says that if a pair of lines is not parallel or orthogonal,
then there are $\omega(n)$ distinct distances between any $n$ points on one line and $n$ points on the other line.
This was proved by Elekes and R\'onyai \cite{ER00}, improved by Elekes \cite{Elek99} to $\Omega(n^{5/4})$ distinct distances, and recently improved to $\Omega(n^{4/3})$ distinct distances by Sharir, Sheffer, and Solymosi \cite{SSS13}.
The approach of \cite{SSS13} was very recently generalized to any two algebraic curves by Pach and De Zeeuw \cite{PdZ13}, who also deduced the following theorem for a single algebraic curve.
This theorem improves the planar case of a recent result by Charalambides \cite{Charal13}, who proved the bound $\Omega_{d,D}(n^{5/4})$ for the number of distinct distances between $n$ points on a common algebraic curve of degree $d$ in $\RR^D$, with an approach based on that of \cite{Elek99}.
It has been conjectured that in all these problems the bound can be improved to $\Omega(n^{2-\varepsilon})$.
\begin{theorem} \label{th:PdZ} {\bf (Pach and De Zeeuw \cite{PdZ13})}
If an algebraic curve in $\RR^2$ of degree $d$ does not contain a line or a circle, then it determines $\Omega_d(n^{4/3})$ distinct distances.
\end{theorem}

In our context, we are interested in the following simple corollary of Theorem \ref{th:PdZ}.
\begin{corollary}\label{co:PdZ}
Let $\pts$ be a set of $n$ points in $\RR^2$ such that $D(\pts)=o(n)$. Then any constant-degree algebraic curve that contains no lines or circles contains $o(n^{3/4})$ points of $\pts$.
\end{corollary}

Thus, our Theorem \ref{th:lines} ``complements" Corollary \ref{co:PdZ} for the case of lines.
Theorem \ref{th:PdZ} cannot be extended to also include lines and circles, since it is possible to place $n$ points on a line or a circle so that they will determine only $\Theta(n)$ distinct distances.
In the proof of Theorem \ref{th:lines} we overcome this problem by showing that either there are many distinct distances between the points that are on the line and the points that are not on the line, or there must be many distinct distances between the points on the line.

To have a bound for every type of curve, it remains to deal with the case of circles, which we do in our second main theorem.
\begin{theorem} \label{th:circles}
Let $\pts$ be a set of $n$ points in $\RR^2$, such that $D(\pts)=o(n)$.
Then there does not exist a circle containing $\Omega(n^{5/6})$ points of $\pts$.
\end{theorem}

In the proof of Theorem \ref{th:circles} we show that there are many distinct distances between the points that are on the circle and the points that are not on the circle.
It follows that if $D(\pts) = o(n)$, then no constant-degree curve contains $\Omega(n^{7/8})$ points of $\pts$. Indeed, we can split such a curve into a constant number of lines and circles, and one curve not containing a line or circle; then we separately apply Theorem \ref{th:lines}, Theorem \ref{th:circles}, and Corollary \ref{co:PdZ} to each of these components.
This leads to the following conjecture, which is related to the question of whether the bound of $\Omega_d(n^{4/3})$ in Theorem \ref{th:PdZ} can be improved to $\Omega_d(n^{2-\varepsilon})$.
\begin{conjecture} \label{con:curves}
Let $\pts$ be a set of $n$ points in $\RR^2$, such that $D(\pts)=o(n)$.
Then there does not exist a constant-degree curve containing $\Omega(n^{1/2+\varepsilon})$ points of $\pts$, for any $\varepsilon>0$.
\end{conjecture}

We conclude the introduction with a few words about the proofs of Theorems \ref{th:lines} and \ref{th:circles}.
Both proofs rely on the approach used in \cite{SSS13} (and subsequently in \cite{PdZ13} and \cite{SS13}), which is based on the Elekes-Sharir framework from \cite{ES11} and \cite{GK11}.
More precisely, this approach defines a set $Q$ of quadruples of points and then applies a double counting argument to $|Q|$. The lower bound for $|Q|$ is obtained by a simple use of the Cauchy-Schwarz inequality.
In \cite{ES11,GK11}, the problem of obtaining an upper bound for $|Q|$ is reduced to a problem about intersections between curves in $\RR^3$, which was solved in \cite{GK11} by introducing a strong new incidence theorem for lines in $\RR^3$. On the other hand, in \cite{PdZ13, SS13, SSS13} the problem of upper bounding $|Q|$ is reduced to a planar incidence problem for curves which can be solved by applying a classic incidence result of Pach-Sharir (see Section \ref{sec:prelims} for more details).
In \cite{PdZ13, SSS13}, this simpler reduction is possible since the points are restricted to lie on a one-dimensional algebraic set.

The above can also be related to the more general theorems of Elekes, R\'onyai, and Szab\'o \cite{Elek02, ER00, ES12}, which say, very roughly, that if a polynomial has fewer values than expected on a Cartesian product of two finite subsets of one-dimensional sets, then the polynomial and the sets must have a special structure.
See also \cite{SS13} for a discussion of this approach.

A new element in our theorems is that only one of the two sets of points is required to lie on a one-dimensional algebraic curve; the other point set is not restricted in any way.
We show that, with more care, the same approach can still be applied in such a case.

To prove Theorem \ref{th:circles}, we use a very recent incidence bound from \cite{WYZ13}, which adapts the Pach-Sharir bound to sets of algebraic curves.
For lines, we also have to deal with many curves that coincide as point sets, for which we use a well-known theorem from additive combinatorics on sumsets and convexity from \cite{ER00}.

In Section \ref{sec:prelims} we introduce the various other results that we will use. In Section \ref{sec:lines} we prove Theorem \ref{th:lines}, and in Section \ref{sec:circles} we prove Theorem \ref{th:circles}.

\section{Preliminaries}\label{sec:prelims}
In this section we introduce some of the tools that we will apply in the proofs of our two main theorems.

Both proofs rely crucially on incidence bounds for algebraic curves.
For a careful definition of algebraic curves we refer to \cite{CLOu}.
Informally, we will think of an algebraic variety as the common zero set in $\RR^m$ or $\CC^m$ of a collection of polynomials. If this zero set has dimension at most one, we will call it an algebraic curve. Frequently, we will deal with algebraic curves in the plane that are the zero set of a single polynomial. If $f\in\RR[z_1,z_2]$ (resp.~$f\in\CC[z_1,z_2]$), we will use $\mathbf{Z}_\RR(f)$ (resp.~$\mathbf{Z}_\CC(f)$) to denote the zero set of $f$. If $f$ is not identically 0, then this will be an algebraic curve (though note that for $f\in\RR[z_1,z_2],$ $\mathbf{Z}_{\RR}(f)$ may be zero-dimensional or contain zero-dimensional components). If $\gamma\subset\RR^2$ is an algebraic curve, we will say that the polynomial $f\in\RR[x,y]$ {\it corresponds} to $\gamma$ if $f$ is the lowest-degree polynomial whose zero set is $\gamma$. We shall do the same for complex curves $\gamma\subset\CC^2.$

The following incidence bound is a well-known result of Pach and Sharir \cite{PS98}, stated for the special case of algebraic curves; see also \cite{KMS12} for a simpler proof of this special case.

\begin{theorem} \label{th:PS} {\bf (Pach and Sharir \cite{PS98})}
Let $\ptss$ be a set of points in $\RR^2$ and let $\curves$ be a set of distinct algebraic curves in $\RR^2$ of degree at most $d$,
such that no two curves of $\curves$ share a one-dimensional component, and any $k$ points of $\RR^2$ are incident to at most $s$ common curves of $\curves$.
Then the number of incidences $I(\ptss,\curves) = |\{(p,\gamma)\in \ptss\times \curves: p\in \gamma \}|$ satisfies
$$I(\ptss,\curves) = O_{d,k,s}\left(|\ptss|^{k/(2k-1)}|\curves|^{(2k-2)/(2k-1)}+|\ptss|+|\curves|\right).$$
\end{theorem}

If a set of curves satisfies the conditions of the theorem, we will say that it has {\it $k$ degrees of freedom} (often neglecting to mention the parameters $d$ and $s$).

In the proof of Theorem \ref{th:lines}, we obtain a multiset $\curves$ of curves that are not necessarily distinct.
We say that the {\it multiplicity} of a curve $\gamma$ is the number of times that $\gamma$ appears in $\curves$.
We will need the following corollary of Theorem \ref{th:PS} for multisets of curves.

\begin{corollary}\label{co:mult}
Let $\ptss$ be a set of points in $\RR^2$ and let $\curves$ be a multiset of algebraic curves with maximum multiplicity $t$,
such that the corresponding set of curves has $k$ degrees of freedom (with parameters $d$ and $s$).
Then
\[ I(\ptss,\curves) =O_{d,k,s}\left(t^{1/(2k-1)}|\ptss|^{k/(2k-1)}|\curves|^{(2k-2)/(2k-1)}+t|\ptss|+\log t|{\curves}|\right).\]
\end{corollary}
\begin{theProof}{\!\!}
For each $0\le i \le \lceil\log_2 t\rceil $, let $\curves_i$ denote the set of curves of $\curves$ that have a multiplicity of at least $2^i$ and at most $2^{i+1}$.
Notice that $|\curves_i|\le |\curves|/2^i$ and that $\sum_{i=1}^{\lceil\log_2 t\rceil}|\curves_i|=|\curves|$.
By Theorem \ref{th:PS}, we have $I(\ptss,\curves_i) = O\big(2^{i+1} (|\ptss|^{\frac{k}{2k-1}}|\curves_i|^{\frac{2k-2}{2k-1}}+|\ptss|+|{\curves_i}|)\big)$. Summing, we obtain
\begin{align*}
I(\ptss,\curves) &=\sum_{i=1}^{\lceil\log_2 t\rceil} O\left(2^{i+1} (|\ptss|^{\frac{k}{2k-1}}(|\curves|/2^i)^{\frac{2k-2}{2k-1}}+|\ptss|+|{\curves_i}|)\right) \\
&=\sum_{i=1}^{\lceil\log_2 t\rceil} O\left(2^{\frac{i+1}{2k-1}} |\ptss|^{\frac{k}{2k-1}}|\curves|^{\frac{2k-2}{2k-1}}+2^{i+1}|\ptss|+2^{i+1}|{\curves_i}|\right) \\
&= O\left(t^{\frac{1}{2k-1}}|\ptss|^{\frac{k}{2k-1}}|\curves|^{\frac{2k-2}{2k-1}}+t|\ptss|+\log t|\curves|\right).
\end{align*}
\end{theProof}


\vspace{20pt}

In the proof of Theorem \ref{th:circles}, we use a variant of Theorem \ref{th:PS} that was recently proved by Wang, Yang, and Zhang \cite{WYZ13}.
It is specific to algebraic curves, and replaces the degrees-of-freedom condition by a more algebraic one, which we now define.

Let $\mathbb{R}\mathbf{P}^N$ be $N$-dimensional real projective space, and let $S_d$ be the set of non-zero polynomials in $\mathbb{R}[x,y]$ of degree at most $d$. We have an injective map $\tau_d\colon S_d\to\RP^N$ given by
\begin{equation}\label{eq:defnOfTauD}
\tau_d\colon \sum a_{ij} x^iy^j \mapsto [a_{ij}].
\end{equation}
If $G\subset S_d$ is a collection of polynomials, we say that $G$ is an {\it algebraic family of dimension $k$} if the image of $G$ under $\tau_d$ is a $k$--dimensional variety. We say that the degree of $G$ is the degree of the corresponding variety. In our applications, we will assume that all varieties of this type are of bounded degree. If $\Gamma$ is a collection of plane algebraic curves, with a corresponding collection of polynomials $G\subset S_d$, we say that $\Gamma$ is {\it contained in an algebraic family of dimension $k$} if $G$ is contained in a bounded-degree algebraic family of dimension $k$.
\begin{theorem} {\bf (Wang, Yang, and Zhang \cite{WYZ13})} \label{th:WYZ}
Let $\ptss$ be a set of points in $\RR^2$ and let $\curves$ be a set of algebraic curves of degree at most $d$. Suppose that $\curves$ is contained in an algebraic family of dimension $k$, and that no two curves in $\curves$ have a common one-dimensional component.
Then
\[ I(\ptss,\curves) = O_{d,k}\left( |\ptss|^{\frac{k}{2k-1}}|\curves|^{\frac{2k-2}{2k-1}}+|\ptss|+|\curves| \right). \]
\end{theorem}

\vspace{20pt}
In the proof of Theorem \ref{th:lines} we will also use a result from additive combinatorics, proved in \cite{ENR00}.
Given sets $A, B \subset \RR$, we set $A+B = \{a+b \mid a\in A, \ b\in B \}$, $A-B = \{a-b \mid a\in A, \ b\in B \}$, and $A^2 = \{a^2 \mid a\in A \}$.

\begin{theorem} [Elekes, Nathanson, and Ruzsa \cite{ENR00}]\label{th:additive}
Let $f:\RR\to\RR$ be a strictly convex or concave function, $A, C, D \subset \RR$, $|A|=n$, and $|C||D|\ge n$. Then we have
\[ |A+C|\cdot|f(A)+D| = \Omega\left(n^{3/2}(|C||D|)^{1/2}\right).\]
\end{theorem}

We note that in \cite{LRN11} this bound was improved slightly, under the condition $|A|\approx |C|$, but this condition is not satisfied in our application.

\section{The proof of Theorem \ref{th:lines}}\label{sec:lines}

Consider a set $\pts$ of $n$ points in $\RR^2$ and a line $\ell$ that contains $\Theta(n^{\alpha})$ points of $\pts$, where $7/8 \le \alpha \le 1$.
We rotate the plane so that $\ell$ is the $x$-axis.
Let $\pts_1 = \pts \cap \ell$, $\pts_2 = \pts \setminus \pts_1$, and let $D= D(\pts_1,\pts_2)$ be the number of distinct distances between $\pts_1$ and $\pts_2$.
We assume, for contradiction, that $D=o(n)$; since $D = D(\pts_1,\pts_2) \le D(\pts)$, obtaining a contradiction to this assumption would complete the proof of Theorem \ref{th:lines}.

For a pair of points $u$ and $v$, we denote by $\|uv\|$ the (Euclidean) length of the straight segment $uv$.
Let $Q$ be the set of quadruples $(a,p,b,q)$ with $a,b \in \pts_1$ and $p,q \in \pts_2$, such that $\|ap\|=\|bq\|$ and $ap \neq bq$.
The quadruples are ordered, so $(a,p,b,q)$ and $(b,q,a,p)$ are considered as two distinct elements of $Q$.

We first derive a lower bound on $|Q|$, as in \cite{SSS13}.
We denote the $D$ distinct distances in $\pts_1\times\pts_2$ as $\delta_1,\ldots, \delta_D$.
Let $E_i = \{(a,p)\in \pts_1 \times \pts_2 \mid \|ap\|=\delta_i \}$, for $i=1,\ldots,D$.
Notice that $\sum_{i=1}^D |E_i| = \Theta(n^{1+\alpha})$.
We have, by the Cauchy-Schwarz inequality and the assumption $D=o(n)$, that
\begin{equation} \label{eq:lowerQbi}
|Q|  = 2\sum_{i=1}^D\binom{|E_i|}{2} \ge \frac{1}{D}\left(\sum_{i=1}^D (|E_i|-1) \right)^2 = \Omega\left(\frac{n^{2+2\alpha}}{D}\right).
\end{equation}

In the remainder of the proof we derive an upper bound on $|Q|$, by reducing the problem to a planar point-curve incidence problem. We partition $Q$ into two parts: $Q^{(1)}$ contains the quadruples $(a,p,b,q)\in Q$ for which $p_y^2=q_y^2$ and $Q^{(2)} = Q \setminus Q^{(1)}$.

We first bound $|Q^{(1)}|$.
In this case, for any choice of the points $a,b,p$, there are at most four choices of $q$ that satisfy both $p_y^2=q_y^2$ and $\|ap\|=\|bq\|$. 
Since there are $O(n^{1+2\alpha})$ possibilities for choosing $a,b,p$, we have $|Q^{(1)}|=O(n^{1+2\alpha})$.

Next, we provide an upper bound for $|Q^{(2)}|$.
A quadruple $(a,p,b,q) \in (\pts_1\times\pts_2)^2$ with $ap\neq bq$ is in $Q$ if and only if $\|ap\|=\|bq\|$, or equivalently,
\begin{equation}
(a_x-p_x)^2+p_y^2 = (b_x-q_x)^2 + q_y^2.
\end{equation}
For any pair of distinct points $p=(p_x,p_y)$ and $q=(q_x,q_y)$ from $\pts_2$,
we define a curve $\gamma_{pq} = \mathbf{Z}_\RR(f_{pq})$ with
\begin{equation}\label{eq:curveDef}
f_{pq} := (x-p_x)^2+p_y^2 - (y-q_x)^2 - q_y^2.
\end{equation}
Since $p_y^2\neq q_y^2$, the curve $\gamma_{pq}$ is a non-degenerate hyperbola.

Let $\ptss$ be the set of all points $(a_x, b_x)$ such that $a=(a_x,0)$ and $b=(b_x,0)$ are points in $\pts_1$.
Let $\curves$ be the multiset of all hyperbolas $\gamma_{pq}$ with $p,q\in \pts_2$.
Then $\gamma_{pq}\in\curves$ is incident to $(a_x,b_x)\in\ptss$ if and only if $(a,p,b,q)\in Q$.
This in turn implies that to obtain an upper bound for $|Q^{(2)}|$, it suffices to obtain an upper bound for the number of incidences between $\curves$ and $\ptss$.
For this, we require the following lemma.
\begin{lemma}
 The multiset of curves $\Gamma$ has maximum multiplicity $O(n^{4(1-\alpha)})$ and $3$ degrees of freedom.
\end{lemma}
\begin{theProof}{\!\!}
We first prove that $\Gamma$ has 3 degrees of freedom. An intersection point of two distinct hyperbolas $\gamma_{pq}$ and $\gamma_{p'q'}$ corresponds to a root of the polynomial $g_{pqp'q'}=f_{pq}-f_{p'q'}$.
By recalling \eqref{eq:curveDef} we notice that $g_{pqp'q'}$ is a linear equation, so the line $\mathbf{Z}(g_{pqp'q'})$ has at most two intersection points with any hyperbola.
This implies that any two distinct hyperbolas $\gamma_{pq}\in\curves$ intersect in at most two points, and thus that any three points are incident to at most one common hyperbola.

We now show that $\Gamma$ has maximum multiplicity  $O(n^{4(1-\alpha)})$. First, note that different pairs of points from $\pts_2^2$ can result in the same hyperbola. The goal is to show that if too many pairs of points result in the same hyperbola, then the underlying point set $\pts$ must contain many distinct distances. By looking at \eqref{eq:curveDef}, we notice that the two pairs $(p,q),(p',q')\in \pts_2^2$ yield the same hyperbola if and only if $p_x=p'_x$, $q_x=q'_x$, and $q_y^2-p_y^2=q'^2_y-p'^2_y$.
For a given pair $(p,q)\in\pts_2^2$, the values of $p'_x,q'_x$ are uniquely determined, and for every choice of $p'_y$ there are at most two valid choices of $q'_y$. Therefore, the maximum multiplicity of a hyperbola is at most twice the maximum number of points of $\pts_2$ on a common vertical line.

Let $\ell_v$ be a vertical line and assume that the set $\pts_v = \pts_2 \cap \ell_v$ contains $\Omega(n^\beta)$ points, where $\alpha +\beta \ge 1$.
We translate the plane so that $\ell \cap \ell_v$ is the origin.
Let $A$ denote the set of $x$-coordinates of points of $\pts_1$ and let $B$ denote the set of $y$-coordinates of points of $\pts_v$ (so $|A|=\Omega(n^\alpha)$ and $|B|=\Omega(n^\beta)$).
Then $A-A$ is the set of distances that are determined by pairs of points of $\pts_1^2$ and $A^2+B^2$ is the set of squares of distances that are determined by pairs of points of $\pts_1\times \pts_v$.
By applying Theorem \ref{th:additive} with $C=-A$, $f(x)=x^2$, and $D=B^2$, we obtain
\[ |A-A|\cdot |A^2 + B^2| = \Omega\left( n^{3\alpha/2}n^{(\alpha+\beta)/2}\right) = \Omega\left( n^{2\alpha+\beta/2}\right). \]
This implies that
$$D(\pts)\ge \max\{D(\pts_1),D(\pts_1,\pts_v)\} =\max\{|A-A|, |A^2+B^2|\} =\Omega\left( n^{\alpha+\beta/4}\right).$$
Since we assumed that $D(\pts)=o(n)$, it follows that no vertical line contains $\Omega(n^{4(1-\alpha)})$ points of $\pts_2$.
That is, the maximum multiplicity of $\curves$ is $O(n^{4(1-\alpha)})$.
\end{theProof}

Applying Corollary \ref{co:mult} with $|\ptss|=O(n^{2\alpha})$, $|\curves|=O(n^2)$, $k=3$, and $t=O(n^{4(1-\alpha)})$ gives
\begin{align*} |Q^{(2)}| &= I(\ptss,\curves) \\
&= O(|\ptss|^{3/5}|\curves|^{4/5}t^{1/5}+\log t|\curves|+t|\ptss|) \\
&= O(n^{6\alpha/5}n^{8/5}n^{4(1-\alpha)/5}+n^2\log n+n^{4-2\alpha})\\
& = O(n^{(2\alpha+12)/5}+n^{4-2\alpha}).
\end{align*}
Combining both cases, we have
\begin{equation*}
\begin{split}
|Q| & = |Q^{(1)}| + |Q^{(2)}| \\
& = O\left(n^{1+2\alpha}+n^{(2\alpha+12)/5}+n^{4-2\alpha}\right).
\end{split}
\end{equation*}
When $\alpha \ge 7/8$, the term $O(n^{1+2\alpha})$ dominates this bound.
Combining it with \eqref{eq:lowerQbi} implies
\[ \frac{n^{2+2\alpha}}{D} = O\left(n^{1+2\alpha}\right), \]
or
\[ D = \Omega\left(n \right).\]
This contradicts the assumption $D=o(n)$ and completes the proof. \smalleop{}
\vspace{2mm}

\noindent {\bf Remark.} It is not difficult to show that $k=3$ is optimal in this case, though the bound $o(n^{4(1-\alpha)})$ on the multiplicity can probably be significantly improved.

\section{The proof of Theorem \ref{th:circles}}\label{sec:circles}

This proof goes along the same lines as the one in Section \ref{sec:lines}.
Consider a set $\pts$ of $n$ points in $\RR^2$ and a circle $C$ that contains $\Theta(n^{\alpha})$ points of $\pts$, where $5/6 \le \alpha \le 1$. We translate and perform a uniform scaling of the plane so that $C$ is the unit circle around the origin.
We also rotate the plane around the origin so that no point of $\pts$ lies on a coordinate axis.
Let $\pts_1 = \pts \cap C$ and $\pts_2 = \pts \setminus \pts_1$,
and let $D = D(\pts_1,\pts_2)$ denote the number of distinct distances between $\pts_1$ and $\pts_2$.
We assume, for contradiction, that $D=o(n)$.

Let $Q$ be the set of quadruples $(a,p,b,q)$ where $a,b \in \pts_1$ and $p,q \in \pts_2$, such that $\|ap\|=\|bq\|$ and $ap \neq bq$.
As in Section \ref{sec:lines}, we have
\begin{align}
|Q| = \Omega\left(\frac{n^{2+2\alpha}}{D}\right). \label{eq:lowerQcirc}
\end{align}
In the remainder of the proof we derive an upper bound on $|Q|$.
We partition $Q$ into two parts: $Q^{(1)}$ contains the quadruples $(a,p,b,q)\in Q$ for which $p_x^2+p_y^2=q_x^2+q_y^2$
(i.e., $p$ and $q$ are on a common circle around the origin), and $Q^{(2)} = Q \setminus Q^{(1)}$.

We first bound $|Q^{(1)}|$. In this case, for any choice of the points $a,b,p$, there are at most two choices of $q$ that satisfy both $p_x^2+p_y^2=q_x^2+q_y^2$ and $\|ap\|=\|bq\|$. Since there are $O(n^{1+2\alpha})$ possibilities for choosing $a,b,p$, we have $|Q^{(1)}|=O(n^{1+2\alpha})$.

Next, we derive an upper bound for $|Q^{(2)}|$, again by reducing the problem to a planar incidence problem.
Unlike in Section \ref{sec:lines}, we first define the curves in four-dimensional space, which makes it easier to prove that no two curves share a component. Then we project the curves to a plane and apply an incidence bound for plane algebraic curves.

A quadruple $(a,p,b,q) \in (\pts_1\times\pts_2)^2$ with $ap\neq bq$ is in $Q$ if and only if $\|ap\|=\|bq\|$, or equivalently,
\begin{equation}\label{eq:condCirc}
(a_x-p_x)^2+(a_y-p_y)^2 = (b_x-q_x)^2 + (b_y-q_y)^2.
\end{equation}
Combining \eqref{eq:condCirc} with the equations $a_x^2+a_y^2=b_x^2+b_y^2=1$ and setting $A_{pq}= (p_x^2+p_y^2-q_x^2-q_y^2)/2$ leads to
\begin{equation*}
p_x a_x+p_y a_y = q_x b_x+q_y b_y + A_{pq}.
\end{equation*}
According to the assumption on $Q^{(2)}$ we have $A_{pq}\neq 0$.

We let $\ptss$ be the set of points $(a,b) = (a_x,a_y,b_x,b_y)\in \CC^4$ such that $a,b\in \pts_1$.
For any pair of points $p,q\in \pts_2$ with $p_x^2+p_y^2\neq q_x^2+q_y^2$, we define a complex curve $\gamma_{pq}\subset \CC^4$ as the set of points $(a_x,a_y,b_x,b_y)\in\CC^4$ satisfying the equations
\begin{equation}
a_xp_x+a_yp_y = b_xq_x+b_yq_y + A_{pq}, \quad\ a_x^2+a_y^2=1, \quad \text{and} \quad b_x^2+b_y^2=1. \label{eq:CurvesCirc}
\end{equation}
Let $\curves$ be the multiset of all curves $\gamma_{pq}$ for which $p,q\in \pts_2$ and $p_x^2+p_y^2\neq q_x^2+q_y^2$.
Notice that $\gamma_{pq}\in \curves$ is incident to $(a,b)\in\ptss$ if and only if $(a,p,b,q)\in Q$. 

We have defined the curves $\gamma_{pq}$ in four-dimensional space because this makes it easier to analyze their intersection properties.
After doing that in Lemma \ref{le:common}, we will project the curves to a plane in Lemma \ref{le:project}, and then apply an incidence bound.
We defined the curves over $\CC$ so that when we project, the resulting sets are also algebraic curves; over $\RR$ this need not be true.



\begin{lemma}\label{le:common}
Any two curves in $\Gamma$ intersect in at most four points.
\end{lemma}
\begin{theProof}{\!\!}
A point $(a_x,a_y,b_x,b_y) \in \gamma_{pq} \cap \gamma_{p'q'}$ satisfies the system of equations
\begin{align}
a_x^2+a_y^2 = 1, &  \qquad \qquad b_x^2+b_y^2 = 1, \label{eq:transform1} \\
p_x  a_x + p_y  a_y &= q_x  b_x + q_y  b_y + A_{pq}, \label{eq:transform2}\\
p_x' a_x + p_y' a_y &= q_x' b_x + q_y' b_y + A_{p'q'}\label{eq:transform3}.
\end{align}
Equations \eqref{eq:transform2} and \eqref{eq:transform3} can be rewritten as
\begin{equation}\label{eq:transformMatrix}
M_{pp'}a
=M_{qq'}b
+A_{pqp'q'},
\end{equation}
where
\begin{equation*}
a =\left[\begin{array}{c}a_x\\a_y\end{array}\right],\qquad b =\left[\begin{array}{c}b_x\\b_y\end{array}\right],\qquad A_{pqp'q'} = [A_{pq},\ A_{p'q'}],
\end{equation*}
\begin{equation*}
M_{pp'}=\left[\begin{array}{ll}p_x & p_y\\ p'_x & p'_y\end{array}\right],\qquad \textrm{and}\qquad M_{qq'}=\left[\begin{array}{ll}q_x & q_y\\ q'_x & q'_y\end{array}\right].
\end{equation*}
We wish to show that there are few points $(a_x,a_y,b_x,b_y)$ that satisfy both \eqref{eq:transform1} and \eqref{eq:transformMatrix}.

First, suppose that $M_{pp'}$ is non-singular. In this case \eqref{eq:transformMatrix} can be rewritten as
\begin{equation}\label{eq:transformMatrixInverted}
a
=M_{pp'}^{-1}M_{qq'}b
+M_{pp'}^{-1}A_{pqp'q'}.
\end{equation}
That is, there is an affine transformation $T_{pqp'q'}$ that takes $b\in C$ to $a\in C$. Since $A_{pq}$ and $A_{p'q'}$ are non-zero, this affine transformation does not fix the origin, and thus the image of $C$ cannot be $C$ itself.
If $T_{pqp'q'}$ is invertible, then $T_{pqp'q'}(C)$ is an algebraic curve of degree 2 different from $C$, so by B\'ezout's theorem (e.g., see \cite[Section 8.7]{CLOu}) we have $|C\cap T_{pqp'q'}(C)|\leq 4$. Then for any $a\in C\cap T_{pqp'q'}(C)$ there is one corresponding $b\in C$, so $|\gamma_{pq}\cap\gamma_{p'q'}|\leq 4$.
If $T_{pqp'q'}$ is singular, then $T_{pqp'q'}(C)$ is contained in a line, hence $|C\cap T_{pqp'q'}(C)|\leq 2$. 
For $a\in C\cap T_{pqp'q'}(C)$, the preimage of $a$ under $T_{pqp'q'}$ is a line, which intersects $C$ in at most 2 points, so there are at most two $b\in C$ that are sent to $a$. Therefore, $|\gamma_{pq}\cap\gamma_{p'q'}|\leq 4$.

If $M_{qq'}$ is non-singular, we can apply the same argument by considering an affine transformation that takes $a\in C$ to $b\in C$.

It remains to consider the case where both $M_{pp'}$ and $M_{qq'}$ are singular. Since it is impossible for more than one of $(p_x,p_y), (p'_x,p'_y), (q_x,q_y), (q'_x,q'_y)$ to be the point $(0,0),$ both $M_{pp'}$ and $M_{qq'}$ must have rank 1. In this case the set $X=\operatorname{range}(M_{pp'}) \cap (\operatorname{range}(M_{qq'}) + A_{pqp'q'})$ is an intersection of two lines, so it is either empty, a single point, or a line.
If $X$ is empty, there are no solutions to \eqref{eq:transformMatrix}, so $|\gamma_{pq}\cap\gamma_{p'q'}|=0$.

Suppose that $X$ consists of a single point $z\in\CC^2$.
We denote by $M_{pp'}^{-1}(z)$ (resp. $M_{qq'}^{-1}(z)$) the set of points $v\in\CC^2$ for which $M_{pp'}v=z$ (resp. $M_{qq'}v=z$). Notice that $M_{pp'}^{-1}(z)$ and $M_{qq'}^{-1}(z)$ are lines in $\CC^2$. Any solution $(a,b) = (a_x,a_y,b_x,b_y)$ to \eqref{eq:transform1} and \eqref{eq:transformMatrix} must satisfy  $a\in C\cap M_{pp'}^{-1}(z)$ and $b\in C\cap M_{qq'}^{-1}(z)$.
This implies that there are at most 2 choices for $a$ and at most 2 choices for $b$, resulting in at most 4 solutions $(a,b)$ to \eqref{eq:transform1} and \eqref{eq:transformMatrix}. Hence, $|\gamma_{pq}\cap\gamma_{p'q'}|\leq 4$.

Finally, suppose that $X$ is a line. This implies that $\operatorname{range}(M_{pp'})$ and $\operatorname{range}(M_{qq'})+A_{pqp'q'}$ are the same line, which is the case if and only if $\operatorname{range}(M_{pp'})$ and $\operatorname{range}(M_{qq'})$ are the same line, and $A_{pqp'q'}$ lies on this line. We will show that then $p=p'$ and $q=q',$ which is forbidden.

Fix a pair of points $p,q$. For another pair $p',q'$ to form the problematic scenario from the previous paragraph with $p,q$, the pair $p',q'$ has to satisfy the following four conditions:
\begin{equation} \label{eq:condits}
\det(M_{pp'})=0,~~\det({M_{qq'}})=0,~~\operatorname{range}(M_{pp'})=\operatorname{range}(M_{qq'}), ~~A_{pqp'q'}\in \operatorname{range}(M_{pp'}).
\end{equation}
The above conditions are equivalent to the equations
\begin{equation*}\label{eq:conditsRePhrased}
\begin{split}
p_xp'_y=p_yp'_x,\qquad
q_xq'_y&=q_yq'_x,\qquad
p_xq'_x=q_xp'_x,\\
p_x(p_x'^2+p_y'^2-q'^2_x-q'^2_y) &= p'_x (p_x^2+p_y^2-q^2_x-q^2_y).
\end{split}
\end{equation*}
We can use the first three equations to eliminate the variables $p_y',q_x',q_y'$ from the fourth, using the fact that $p_x,p_x'\neq 0$ (since no point of $\pts$ lies on an axis), and also $A_{pq}, A_{p'q'}\neq 0$.
This leads to $p_x' = p_x$, which then implies $p_y=p_y',q_x=q_x'$, and $q_y=q_y'$.
%
%
%
%
%
\end{theProof}


\begin{lemma}\label{le:project}
There exist a point set $\ptss^*\subset \RR^2$ and a set $\curves^*$ of algebraic curves in $\RR^2$ of degree at most 4, with $|\ptss^*|=|\ptss|=O(n^{2\alpha})$ and $|\curves^*|=|\curves|=O(n^2)$, such that:
\begin{itemize}
\item $\curves^*$ is contained in an algebraic family of dimension 4;
\item No two curves have a common one-dimensional component;
\item Every quadruple in $Q^{(2)}$ corresponds to an incidence between a point in $\ptss^*$ and a curve in $\curves^*$,
and every such incidence corresponds to at most four quadruples.
\end{itemize}
\end{lemma}
\begin{theProof}{\!\!}
Let $\mathcal O$ be a generic rotation in $\RR^4$ (i.e.~a generic element of the orthogonal group $O(4,\RR)$). Note that $\mathcal{O}$ can also be identified with an element $\mathcal{O}^{\CC}$ of $O(4,\CC)$ via the standard embedding of $O(4,\RR)$ into $O(4,\CC).$ Let $\pi\colon\RR^4\to\RR^2$ be the projection $(x_1,...,x_4)\mapsto(x_1,x_2)$. By abuse of notation, we will also use $\pi$ to refer to the corresponding projection from $\CC^4$ to $\CC^2$.
For each $\gamma_{pq}\in\curves$, the set $\pi(\mathcal{O}^{\CC}(\gamma_{pq}))$ is a plane algebraic curve of degree at most 4.
We define $\ptss^* = \pi(\mathcal{O}(\ptss))$ and note that $\ptss^* \subset \RR^2$ since $\ptss \subset \RR^4$.

Let $f_{pq}(z_1,z_2)\in\CC[z_1,z_2]$ be the polynomial corresponding to $\pi(\mathcal{O}^{\CC}(\gamma_{pq}))$. Then $f_{pq}$ is a polynomial of degree at most 4, and $\mathbf{Z}_\CC(f_{pq})=\pi(\mathcal{O}^{\CC}(\gamma_{pq}))$. Note that the polynomial $f_{pq}$ has real coefficients, since it can be obtained by taking a bivariate resultant of three polynomials with real coefficients (for an introduction to resultants see \cite[Section 3.5-3.6]{CLOu}). These three polynomials are  the three polynomials from \eqref{eq:CurvesCirc} after they have been pre-composed with $\mathcal{O}.$ Since no two curves $\gamma_{pq},$ $\gamma_{p'q'}\in\Gamma$ share a component, and since $\mathcal{O}$ is generic, no two curves $\mathbf{Z}_\CC(f_{pq})$ and $\mathbf{Z}_\CC(f_{p'q'})$ share a component\footnote{Note that $\mathcal{O}^{\CC}$ is the complexification of a generic element of generic element of $O(4,\RR),$ which doesn't imply that $\mathcal{O}^{\CC}$ is a generic element of $O(4,\CC)$. However, it does imply 
that the projections $\pi(\mathcal{O}^{\CC}(\gamma_{pq}))$ and $\pi(\mathcal{O}^{\CC}(\gamma_{p'q'}))$ do not share a component.}. Also note that since $\pi\circ \mathcal{O}$ is generic, we have that $|\pi^{-1}(z)\cap \gamma_{pq}|\leq 4$ for any $z\in\CC^2$.

We now consider $\mathbf{Z}_\RR(f_{pq})$, which is the set of real points of $\mathbf{Z}_\CC(f_{pq})$. If $\mathbf{Z}_\RR(f_{pq})$ were all of $\RR^2$, then we must have $\mathbf{Z}_\CC(f_{pq})=\CC^2,$ which we know cannot be the case, since we can verify that the three varieties in $\CC^4$ defined by the three polynomials from \eqref{eq:CurvesCirc} intersect properly. Thus $\mathbf{Z}_\RR(f_{pq})\subset\RR^2$ must be a plane algebraic curve of degree at most 4.
We set
$$\Gamma^* = \{\mathbf{Z}_\RR(f_{pq}): \gamma_{pq}\in \Gamma\}.$$
Note that no two curves $\mathbf{Z}_\RR(f_{pq})$ and $\mathbf{Z}_\RR(f_{pq})$ can share a one-dimensional component because no two curves $\mathbf{Z}_\CC(f_{pq})$ and $\mathbf{Z}_\CC(f_{p'q'})$ share a component. Since $|\pi^{-1}(z)\cap \gamma_{pq}|\leq 4$ for any $z\in\CC^2$, at most four quadruples correspond to the same incidence of $\mathbf{Z}_\RR(f_{pq})$ with $\ptss^*$.

Finally, we show that $\Gamma^*$ is contained in an algebraic family of dimension 4.
Note that the polynomial $f_{pq}$ is of the form
\begin{equation}\label{eq:formOfFpq}
f_{pq}(z_1,z_2) = \sum_{\substack{i,j\geq 0\\ i+j\leq 4}}a_{ij}(p_x,p_y,q_x,q_y)z_1^iz_2^j,
\end{equation}
where each function $a_{ij}$ is a bounded degree polynomial in $p_x,p_y,q_x,$ and $q_y$ (the functions $a_{ij}$ depend on our choice of $\mathcal{O}$).
Recall the map $\tau_4\colon S_4\to \RP^{14}$ defined in \eqref{eq:defnOfTauD}.
Let $S'\subset S_4$ be the set of all polynomials of the form \eqref{eq:formOfFpq} for all $(p_x,p_y,q_x,q_y)\in\RR^4$.  Then
\begin{equation}
\tau_4(S') = \{ [a_{ij}(p_x,p_y,q_x,q_y)]_{\substack{i,j\geq 0\\ i+j\leq 4}}\colon (p_x,p_y,q_x,q_y)\in\RR^4 \}.
\end{equation}
Thus $\tau_4(S')\subset\RP^{14}$ is a 4--dimensional variety of bounded degree, so $S'$ is an algebraic family of dimension 4.
Each curve in $\Gamma^*$ corresponds to a polynomial in $S'$, so $\Gamma^*$ is contained in an algebraic family of dimension 4.
\end{theProof}
\vspace{2mm}

We may now apply Theorem \ref{th:WYZ} with $k=4$ to obtain the bound
\begin{align*}
|Q^{(2)}| &\leq  4 I(\ptss^*, \curves^*) \\
&=O\left(|\ptss|^{4/7}|\curves|^{6/7}+|\ptss|+|{\curves}|\right) \\
&= O\left(n^{8\alpha/7}n^{12/7}+n^{2\alpha}+n^2\right) \\
&= O\left(n^{(8\alpha+12)/7} \right).
\end{align*}
Combining both cases, we have
\begin{equation*}
\begin{split}
|Q| & = |Q^{(1)}| + |Q^{(2)}| \\
& = O\left(n^{1+2\alpha}+n^{(8\alpha+12)/7}\right).
\end{split}
\end{equation*}
When $\alpha \ge 5/6$, the term $O(n^{1+2\alpha})$ dominates this bound.
As in Section \ref{sec:lines}, combining this bound with \eqref{eq:lowerQcirc} gives $D = \Omega(n)$, which contradicts the assumption $D = o(n)$ and completes the proof of Theorem \ref{th:circles}.

\vspace{20pt}
\noindent {\bf Remarks.} It is not difficult to show that in this case $k=4$ cannot be improved on.
Note that Lemma \ref{le:common} already implied that the curves have five degrees of freedom, which would lead to a weaker bound using $k=5$ in Theorem \ref{th:PS}.
On the other hand, it should be possible to use Theorem \ref{th:WYZ} to also prove Theorem \ref{th:lines}, but there it was considerably easier to prove directly that the curves had three degrees of freedom.

\vspace{20pt}\noindent
{\bf Acknowledgement}\\
The third author would like to thank J\'anos Pach and Filip Mori\'c for helpful discussions regarding this problem.


\end{document}